\documentclass[12pt]{article}
\usepackage[english]{babel}
\usepackage[T1]{fontenc}
\usepackage[latin1]{inputenc}
\usepackage{fullpage}
\usepackage{amsfonts}
\usepackage[all]{xy}
\usepackage{stmaryrd}
\usepackage{xypic}
\usepackage{psfrag}
\usepackage{epsfig}

\usepackage{amsmath}
\usepackage{amssymb,latexsym}
\usepackage[mathscr]{eucal}
\usepackage{url}
\usepackage{fancybox}

\newtheorem{thm} {\bf  Theorem\bf} [subsection]
\newtheorem{cor} [thm] {\bf  Corollary\bf}
\newtheorem{defi} [thm] {\bf  Definition\bf}

\newtheorem{prop} [thm] {\bf  Proposition\bf}

\newtheorem{rem}[thm]{\bf  Remark\bf}
\newtheorem{lem}[thm]{\bf  Lemma\bf}

 \newcommand{\OO}{\mathcal{O}}
 \newcommand{\su}{\mathcal{SU}_C(2)}

 \newcommand{\ra}{\rightarrow}

 \newcommand{\til}{\tilde{\phi}}
 \newcommand{\ZZ}{\mathbb{Z}}

\newcommand{\lra}{\longrightarrow}
\newcommand{\pr}{\mathbb{P}}

\newcommand{\A}{\mathcal{A}}

\newcommand{\pro}{\pr^4_{\omega}}
\newcommand{\pru}{\pr^3_{\omega+}}
\newcommand{\prd}{\pr^2_{\omega}}

\begin{document}

\begin{center}
\LARGE{A conic bundle degenerating on the Kummer surface}
\end{center}

\vspace{2pt}

\begin{center}
\large{\textsc{Michele Bolognesi}}
\end{center}
\begin{abstract}
\vspace{2pt} \footnotesize{ \noindent Let $C$ be a genus 2 curve and $\su$ the moduli space of semi-stable rank
2 vector bundles on $C$ with trivial determinant. In \cite{bol:wed} we described the parameter space of non
stable extension classes (invariant with respect to the hyperelliptic involution) of the canonical sheaf
$\omega$ of $C$ with $\omega_C^{-1}$. In this paper we study the classifying rational map $\varphi: \pr
Ext^1(\omega,\omega^{-1})\cong \pr^4 \dashrightarrow \su\cong \pr^3$ that sends an extension class on the
corresponding rank two vector bundle. Moreover we prove that, if we blow up $\pr^4$ along a certain cubic
surface $S$ and $\su$ at the point $p$ corresponding to the bundle $\OO \oplus \OO$, then the induced morphism
$\tilde{\varphi}: Bl_S \ra Bl_p\su$ defines a conic bundle that degenerates on the blow up (at $p$) of the
Kummer surface naturally contained in $\su$. Furthermore we construct the $\pr^2$-bundle that contains the conic
bundle and we discuss the stability and deformations of one of its components.

}


\end{abstract}

\parindent=0pt













\section*{Introduction}

Let $C$ be a smooth genus 2 curve. Let $Pic^1(C)$ be the Picard variety that parametrizes all degree 1 line
bundles on $C$ and $\Theta$ the canonical theta divisor made up set-theoretically of line bundles $L$ s.t.
$h^0(C,L)\neq 0$. We will denote $\su$ the moduli space of semi-stable rank 2 vector bundles on $C$ with trivial
determinant. The description of this moduli space dates back almost fourty years ago to the paper
\cite{rana:cra}. Ramanan and Narashiman proved that $\su$ is isomorphic to the linear system $|2\Theta|\cong
\pr^3$ on $Pic^1(C)$ and that the semi-stable locus is exactly the Kummer quartic surface image of the Jacobian
$Jac(C)$ via the Kummer map. In this paper we look at $\su$ in a different frame. Let $\omega$ be the canonical
bundle on $C$, we consider the space $\pr Ext^1(\omega, \omega^{-1})=: \pro = |\omega^3|^*$. This space
parametrizes extension classes $(e)$ of $\omega$ by $\omega^{-1}$.

$$0 \lra \omega \lra E_e \lra \omega^{-1} \lra 0\ \ \ \ \ \ \ (e)$$

Therefore there exists a classifying map

$$\varphi:\pro \dashrightarrow \su$$

that associates the vector bundle $E_e$ to the extension class $(e)$. Bertram showed in \cite{ab:rk2} that
$\varphi$ is given by the quadrics in the ideal $\mathcal{I}_C(2)$ of the curve, that is naturally embedded as a
sextic in $\pro$. In the first part of the paper we describe the fibers of the map $\varphi$. Let $E\in\su$ and
$C_E$ the closure $\overline{\varphi^{-1}(E)}$ of the fiber of $E$, then the principal results of section 1 are
the following.

\begin{thm}
Let $E$ be a stable vector bundle, then $dim C_E=1$ and $C_E$ is a smooth conic.

\end{thm}

\begin{thm}
Let $E$ be a strictly semi-stable vector bundle, then $C_E$ is singular. If $E\cong L\oplus L^{-1}$ for $L\in
JC[2]/\mathcal{O}_C$ then $C_E$ is a double line.
\end{thm}

Moreover we show that the fiber over the $S$-equivalence class of the bundle $\OO_C \oplus \OO_C$ is a cone
$S\in \pro$ over a twisted cubic curve. In Section 2 we blow up $\pro$ along the surface $S$ and $|2\Theta|$ at
the origin of the Kummer surface $K^0$ that represents the semi-stable boundary. Let $Bl_S(\pro)$ be the blow up
of $\pro$ and $\pr^3_0$ the blow-up of $|2\Theta|$, we describe the induced map $\tilde{\varphi}:Bl_S(\pr^4)\lra
\pr^3_{\OO}$ and we prove that the restriction of $\tilde{\varphi}$ to the exceptional divisors is a conic
bundle. The main theorem of Section 2 is in fact the following.

\begin{thm}
The morphism

$$\tilde{\varphi}:Bl_S(\pr^4)\lra \pr^3_{\OO}$$

is a conic bundle whose discriminant locus is the blow-up at the origin of the Kummer surface $K^0$.
\end{thm}

Moreover we construct a rank 2 vector bundle $\mathcal{A}$ on $\pr^2$, such that $\pr (\mathcal{A}\oplus
\OO_{\prd}) \cong Bl_S(\pr^4)$. This leads us to prove that our conic bundle can be seen as a section

$$\tilde{\varphi}\in  \mathrm{Hom}(\OO_{\prd}\oplus\OO_{\prd}(-1),Sym^2 (\mathcal{A}^*\oplus \OO_{\prd})
\otimes \OO_{\prd}(-1)).$$

This vector space has dimension 16 thus we have a moduli map that associates to a smooth genus 2 curve the conic
bundle induced by its classifying map.

\begin{eqnarray*}
\Xi:\{ \mathrm{smooth\ genus\ 2\ curves} \} & \lra & \pr^{15}=\pr B;\\
C & \mapsto & \tilde{\varphi}_C.
\end{eqnarray*}

In Section 3 we study the stability and the deformations of the bundle $\mathcal{A}$. By applying the Hoppe
criterion (Prop. \ref{kit}), we get the following Theorem.

\begin{thm}
The vector bundle $\A$ on $\prd$ is stable.
\end{thm}

Finally, via a few cohomology calculations we find the dimension of $Ext^1(\mathcal{A},\mathcal{A})$.

\begin{thm}
The space of deformations of $\A$ has dimension $dim(Ext^1(\A,\A))=5$.
\end{thm}

\textit{Acknowledgments:} I would like to thank my Phd advisor Christian Pauly for his suggestions, Chiara
Brambilla for a useful advice and Yves Laszlo for posing the problem of the deformations of the vector bundle
$\mathcal{A}$.

\section{The classifying map}

\subsection{Preliminaries on extension classes}

Let $C$ be a smooth genus 2 curve and $\lambda$ the hyperelliptic involution on $C$; we will denote
$\mathcal{W}$ the set of the Weierstrass points of $C$. Let also $Pic^d(C)$ be the Picard variety parametrizing
degree $d$ line bundles over $C$ and $Jac(C)=Pic^0(C)$ the Jacobian variety of $C$. We will denote $K^0$ the
Kummer surface obtained as quotient of $Jac(C)$ by $\pm Id$ and $K^1$ the quotient of $Pic^1(C)$ by the
involution $\tau:\xi\mapsto \omega \otimes \xi^{-1}$. Moreover we remark that the 16 theta characteristics are
the fixed points of the involution $\tau$. Let $\Theta\subset Pic^1(C)$ be the Riemann theta divisor. It is
isomorphic to the curve $C$ via the Abel-Jacobi embedding

\begin{eqnarray}\label{eq:aj}
Aj: C & \hookrightarrow & Pic^1(C), \\
p & \mapsto & \OO_C(p). \nonumber
\end{eqnarray}

Let $\su$ be the moduli space of semi-stable rank two vector bundles on $C$ with trivial determinant. It is
isomorphic to $\pr^3\cong |2\Theta|$, the isomorphism being given by the map \cite{bo:fib1}

\begin{eqnarray}\label{tetta}
\theta: \su & \lra & |2\Theta|,\\
E & \mapsto & \theta(E);\nonumber
\end{eqnarray}

where

$$\theta(E):=\{L \in Pic^1(C)| h^0(C,E\otimes L)\neq 0\}.$$

\bigskip

With its natural scheme structure, $\theta(E)$ is in fact linearly equivalent to $2\Theta$. The Kummer surface
$K^0$ is embedded in $|2\Theta|$ and points in $K^0$ correspond to bundles $E$ whose S-equivalence class $[E]$
contains a decomposable bundle of the form $M\oplus M^{-1}$, for $M\in Jac(C)$. Furthermore on the semistable
boundary the morphism $\theta$ restricts to the Kummer map.

Let $\omega$ be the canonical line bundle on $C$. We introduce the 4-dimensional projective space

$$\pr^4_{\omega}:=\pr
Ext^1(\omega,\omega^{-1})=|\omega^3|^*.$$

A point $e \in \pro$ corresponds to an isomorphism class of extensions

$$0\lra \omega^{-1} \lra E_e \lra \omega\lra 0.\ \ \ \ \ \ (e)$$

We denote by $\varphi$ the rational classifying map

\begin{eqnarray*}
\varphi:\pro & \dashrightarrow & |2\Theta|\\
e & \mapsto & \textrm{S-equivalence class of } E_e.
\end{eqnarray*}

Let $\mathcal{I}_C$ be the ideal sheaf of the curve $C\subset \pro$, Bertram (\cite{ab:rk2}, Theorem 2) showed
that there is an isomorphism (induced via pull-back by $\varphi$)

$$H^0(\su,\mathcal{O}(2\Theta))\cong H^0(\pro,\mathcal{I}_C\otimes\mathcal{O}(2)).$$

Therefore the classifying map $\varphi$ is the rational map given by the full linear system of quadrics
contained in the ideal of $C \subset \pro$. In fact the locus of non semi-stable extensions is exactly
represented by $C$, as the next lemma shows.

\begin{lem}\label{lem:ber}\cite{ab:rk2}
Let $(e)$ be an extension class in $\pro$ and $Sec(C)$ the secant variety of $C\subset \pro$, then the vector
bundle $E_e$ is not semistable if and only if $e \in C$ and it is not stable if and only if $e \in Sec(C)$.
\end{lem}

\begin{rem}\label{re:marco}
One can say even more. In fact, given $x,y\in C$ the secant line $\overline{xy}$ is the fiber of $\varphi$ over
the S-equivalence class of
$\omega(-x-y)\oplus\omega^{-1}(x+y)$.\\

\end{rem}

This implies directly the following Corollary.

\begin{cor}

The image of the secant variety Sec(C) by the classifying map $\varphi$ is the Kummer surface $K^0\subset
|2\Theta|$.

\end{cor}

The hyperelliptic involution $\lambda$ acts on the canonical line bundle over $C$ and on its spaces of sections.
A straightforward Riemann-Roch computation shows that $h^0(C,\omega^3)^*=5$. Let $\pi: C\rightarrow \pr^1$ be
the hyperelliptic map. There is a canonical linearization for the action of $\lambda$ on $\omega$ that comes
from the fact that $\omega=\pi^*\mathcal{O}_{\pr^1}(1)$. In fact, by Kempf's Theorem (\cite{dn:pfv},
Th\'{e}or\`{e}me 2.3), a line bundle on $C$ descends to $\pr^1$ if and only if the involution acts trivially on
the fibers over Weierstrass points. Thus we choose the linearization
$\delta:\lambda^*\omega\stackrel{\sim}{\rightarrow} \omega$ that induces the identity on the fibers over
Weierstrass points. This means that

$$Tr(\lambda: L_{w_i}\rightarrow L_{w_i})=1,$$

for every Weierstrass point $w_i$. Moreover we have that $d\lambda_{w_i}=-1$, which implies, via the
Atiyah-Bott-Lefschetz fixed point formula (\cite{gh:pag}, p.421), that

$$h^0(C,\omega^3)_+ - h^0(C,\omega^3)_-=3.$$

Since $ h^0(C,\omega^3)_+ + h^0(C,\omega^3)_-=5$, this means that $h^0(C,\omega^3)_+=4$ and
$h^0(C,\omega^3)_-=1$ and we can see that

$$H^0(C,\omega^3)_-=\sum_{i=1}^6 w_i.$$

Furthermore, we have

$$E_{\lambda(e)}=\lambda^*E_e$$

thus the points of $\pru:= \pr H^0(C,\omega^3)_+^*$ represent involution invariant extension classes. We have
studied this classes in \cite{bol:wed}; in this paper our aim is to describe more precisely the classifying map.

\subsection{The fibers of the classifying map}

Let $p\in \pr^3\cong |2\Theta|$ be a general point, then the fiber $\varphi^{-1}(p)$ is the intersection of 3
quadrics

$$Q_1\cap Q_2 \cap Q_3 = \varphi^{-1}(p).$$

If $dim (\varphi^{-1}(p))=1$ then the fiber $\varphi^{-1}(p)$ is a degree 8 curve. Since $C\subset
\varphi^{-1}(p)$
and $deg(C)=6$, the residual curve $\varphi^{-1}(p)- C$ is a conic.\\

We will often denote $V$ the space of global sections $H^0(C,\omega).$

We have the following equality

$$H^0(C,\omega^3)_+=Sym^3V$$

and we denote $X\subset \pr (Sym^3 V^*)$ the twisted cubic curve image of $\pr V \cong \pr^1$ via the cubic
Veronese morphism. Let's consider now the linear subspace $<D>$ of $\pro$ generated by the points of a divisor
$D\in |\omega^2|$. Since $D\in |\omega^2|$, we can write it down as

$$D= a + b + \lambda(a) + \lambda (b)$$

for $a,b \in C$. Furthermore we remark that the annihilator of $<D>$ is
$H^0(C,\omega^3(-a-b-\lambda(a)-\lambda(b)))$, that has dimension equal to 2. This means that the linear envelop
$<D>$ is a $\pr^2$, and we shall denote it as $\pr^2_{ab}$.

\begin{prop}\cite{ln85}\label{carol}
Let $c,d \in C$ and $e \in \pro$ an extension

$$0\lra \omega^{-1}\stackrel{i_e}{\lra} E_e \stackrel{\pi_e}{\lra}\omega\lra 0.$$

Then $e \in \pr^2_{cd}$ if and only if it exists a section $\beta \in H^0(C,Hom(\omega^{-1},E))$ s.t.

$$Zeros(\pi_e \circ \beta)=c+d+ \lambda(c) + \lambda (d).$$

\end{prop}

We will denote $[\OO_C \oplus \OO_C]$ the $S$-equivalence class of the rank 2 bundle $\OO_C \oplus \OO_C$.

\begin{prop}\label{pr:giu}

Let

$$|\omega^3|^*=\pro \stackrel{\varphi}{\dashrightarrow}\pr^3=|2\Theta|$$

be the classifying map. Then the closure of $\varphi^{-1}([\OO_C \oplus \OO_C])$ is the cone $S$ over a twisted
cubic curve $X\subset \pru$.

\end{prop}

\textit{Proof:} The vertex of $S$ is the point $x=\pr H^0(C,\omega^3)^*_-\in \pro$, that is the projectivized
anti-invariant eigen-space. This means that every line contained in $S$ is a secant of $C$ invariant under the
involution of $\pro$. Such a secant line can be written as $\overline{p\lambda(p)}$. The image of such a secant
line via $\varphi$ is the origin, hence $S \subset \varphi^{-1}([\OO_C \oplus \OO_C])$. In order to prove the
opposite inclusion we remark that a vector bundle $E$ contained in the $S$-equivalence class of the origin
satisfies the following exact sequence

$$0\lra \OO_C\stackrel{\nu_E}{\lra} E \lra \OO_C \lra 0.$$

This implies that the trivial bundle is a sub-bundle of $E$.

Let us consider the morphism $\varsigma$, composition of $\nu_E$ and $\pi_E$:

$$\varsigma: \OO_C \stackrel{\nu_E}{\lra} E \stackrel{\pi_E}{\lra} \omega.$$

Then $\delta \in H^0(C,\omega)=Hom (\OO_C,\omega)$ and we obtain the following diagram.

$$\begin{array}{ccccccccc}
& & & & \OO_C & & & &\\
& & & & \downarrow &\stackrel{\searrow}{\varsigma} & & &\\
0 & \ra & \omega^{-1} & \ra & E & \ra & \omega & \ra & 0.\\
\end{array}$$

This means that the morphism $\varsigma$ is a section of $H^0(C,\omega)$ and its divisor is of type $a +
\lambda(a)$ for $a \in C$. Because of Proposition \ref{carol} this implies that the extension classes contained
in the fiber $\varphi^{-1}([\OO_C \oplus \OO_C])$ belong to an invariant secant line, therefore
$\varphi^{-1}([\OO_C \oplus \OO_C])\subset S \square$\\

We will see that the map $\varphi$ defines a conic bundle on $|2\Theta|- [\OO_C \oplus \OO_C]$. We will now
describe the
fibers of $\varphi$ on the open set complementary to the origin.\\

Let $E\in\su$ we will call $C_E$ the closure $\overline{\varphi^{-1}(E)}$ of the fiber of $E$.

\begin{thm}\label{yo}
Let $E$ be a stable vector bundle, then $dim C_E=1$ and $C_E$ is a smooth conic.

\end{thm}

\textit{Proof:} Let $E$ be a stable bundle and $e$ the following equivalence class of extensions

\begin{equation}\label{eq:fly}
0\lra\omega^{-1}\lra E_e\cong E \lra \omega \lra 0
\end{equation}

For a general bundle $E$, $C_E$ has dimension equal to 1. In fact for a genus 2 curve the Riemann-Roch theorem
gives $\chi(E\omega)=4+2(-1)=2$. Moreover, by Serre duality, we have $h^1(E\omega)=h^0(E^*)$ and
$h^0(E^*)=h^0(Hom(E,\mathcal{O}_C))=0$ because $E$ is stable.\\

We define the following map

\begin{eqnarray*}
j:\varphi^{-1}(E) & \lra & \pr H^0(C,E\omega)=\pr^1,\\
 e & \mapsto & j(e),
\end{eqnarray*}

that sends the extension class $e\in C_E$ on the point of $\pr H^0(C,E\omega)$ corresponding to the first morphism of the
exact sequence \ref{eq:fly}. This map has degree 1 and it is not defined on the points of $C_E\cap C$.\\

Furthermore we remark that the projection from $\su$ with centre $[\OO_C \oplus \OO_C]$ can be described in the
following way.

\begin{eqnarray*}
\Delta:\su & \dashrightarrow & \pr^2 = |\omega^2|,\\
E & \mapsto & D(E).\\
\end{eqnarray*}

Here by $D(E)$ we mean the divisor on $C$ with the following support.

$$Supp(D(E))=\{p\in C| h^0(E\otimes \mathcal{O}_C(p))\neq 0\}.$$

In fact the projection from $\OO$ is exactly the restriction to $C$, embedded in $Pic^1(C)$ via the Abel-Jacobi
 map of equation \ref{eq:aj}, of the map $\theta$ from equation \ref{tetta}. Now we consider the determinant map

\begin{eqnarray*}
\bigwedge^2 H^0(E\omega) & \lra & H^0(\omega^2),\\
s\wedge t & \mapsto & Zero(s\wedge t).\\
\end{eqnarray*}

Let $p\in C$ be a point of the curve, if $p\in$ Zero$(s \wedge t)$ then there exists a non zero section $s_p\in
H^0(C,E\omega (-p))$. Hence $h^0(C,E\omega(-p))\neq 0$. Moreover, if we make the hyperelliptic involution
$\lambda$ act on $E\omega(-p)$ we find that $h^0(C,E\omega(-p))=h^0(C,E\otimes \OO_C(p))\neq 0$. This implies
that the zero divisor of $s\wedge t$ is $D(E)$. Now $D(E)$ has degree 4 and for every $p\in D(E)$ it exists a
section $s_p$. We remark then that the morphism $j$ is surjective on the open set $\pr H^0(C,E\omega)/
\{s_p|p\in D(E) \}$, that means that it is
dominant.\\

In order to end the proof we need 3 technical lemmas.

\begin{lem}\label{shit}
Let $D\in |\omega^2|$ be the divisor $a+ b + \lambda(a) + lambda (b)$. The image of
$$\varphi_{|\pr^2_{ab}}:\pr^2_{ab} \dashrightarrow \pr^3$$
is the fiber of $\Delta$ over $D$, i.e. the line passing by $[\OO_C\oplus \OO_C]$ and the point corresponding to
$D$ in $|\omega^2|$.
\end{lem}

\textit{Proof:} We remark first that $D=C\cap \pr^2_{ab}$, because $H^0(C,\omega^3(-a-b-\lambda(a)-\lambda(b)-
c))=H^0(C,\omega - c)=1$ for every $c\in C$. This implies that the restriction

$$\varphi_{|\pr^2_{ab}}:\pr^2_{ab} \dashrightarrow \pr^3$$

is given by quadrics passing by the four points of $D$. The space of quadrics on $\pr^2$ has dimension 5 and we
impose 4 independent linear conditions. So the image of $\pr^2_{ab}$ via $\varphi_{|\pr^2_{ab}}$ is a
$\pr^1\subset \pr^3$.

Let $e$ be an extension class in $\pr^2_{ab}$ and $E_e$ its image via $\varphi$ in $\su$. Now, by Proposition
\ref{carol}, the extension $e$ belongs to $\pr^2_{ab}$ if and only if it exists a section $\alpha \in
H^0(C,Hom(\omega^{-1},E))$ s.t., using the notation of the following diagram, we have Zeros$(\pi_e \circ
\alpha)=D$.

$$\begin{array}{ccccccccc}
 & & & & \omega^{-1} & & & &\\
 & & & & \alpha \downarrow & \searrow & & & \\
 0 & \ra & \omega^{-1} & \stackrel{i_e}{\ra} & E_e & \stackrel{\pi_e}{\ra} & \omega & \ra 0. \\
\end{array}$$

This implies that $\alpha$ and $i_e$ are 2 independent sections of $E\omega$ and Zeros$(i_e\wedge \alpha)= D$. $\square$\\

\begin{lem}\label{anto}

Let $E\in \su$, then we have the equality

$$C_E\cap C= D(E).$$
\end{lem}

\textit{Proof:} Let $c,d \in C$ and let us suppose that $D(E)=c+d+\lambda(c)+\lambda(d)$.  We remind that
$\pr^2_{cd}$ is the plane s.t. $c+d+\lambda(c)+\lambda(d)\subset \pr^2_{cd}$. By Lemma \ref{shit} the fiber
$C_E$ is a conic contained in $\pr^2_{cd}$ and passing by the 4 points of $D(E)$. Since $D(E)=C\cap \pr^2_{cd}$
we obtain the equality $C_E\cap C= D(E)$.$\square$

\begin{lem}\label{ss}

Let $E\in \su$ be a stable bundle, then we have a decomposition

$$H^0(C,E\omega)=H^0(C,E\omega)_+\oplus H^0(C,E\omega)_-$$

of $H^0(C,E\omega)$ in two eigen-spaces of dimension 1.
\end{lem}

\textit{Proof:} We will use the Atiyah-Bott-Lefschetz formula, so we must choose a linearization

$$\nu:\lambda^*E \stackrel{\sim}{\lra} E$$

for the action of $\lambda$ on $E$ and watch how it acts on the fibers over the points of $\mathcal{W}$. The
bundle $E$ is an extension

$$0\lra \omega^{-1} \lra E \lra \omega \lra 0$$

so a linearization on $E$ is defined once one chooses two linearizations on $\omega $ and $\omega^{-1}$. We have
already chosen

$$\delta:\lambda^*\omega \lra \omega$$

that acts trivially on the fibers over the points of $\mathcal{W}$. For the line bundle $\omega^{-1}$ we have
two different choices: the linearization that acts trivially on the fibers over the points of $\mathcal{W}$ and
its inverse. Let $x\in \mathcal{W}$, then we can decompose

\begin{equation}\label{deco}
E_x=\omega_x \oplus \omega^{-1}_x.
\end{equation}

If we choose the first linearization on $\omega^{-1}$  then, by Kempf Lemma, the vector bundle $E$ would be the
pull-back of a bundle $F$ over $\pr^1$ defined by the exact sequence

\begin{equation}\label{pale}
0\lra \OO_{\pr^1}(-1) \lra F \lra \OO_{\pr^1}(1) \lra 0.
\end{equation}

We remark that the only vector bundle on $\pr^1$ that verifies the exact sequence (\ref{pale}) is
$\OO_{\pr^1}(-1)\oplus \OO_{\pr^1}(1)$ and it is not semi-stable. Then the choice of $\delta$ as a linearization
for $\omega$ forces us to choose the linearization

$$\tilde{\delta}:\lambda^*\omega^{-1} \stackrel{\sim}{\lra} \omega^{-1}$$

that induces $-Id$ on the fibers over the points of $\mathcal{W}$. We recall that $h^0(C,E\omega)=2$. Thanks to
the decomposition (\ref{deco}) the trace of the linearization $\delta \oplus \tilde{\delta}$ is zero. Then by
the Atiyah-Bott-Lefschetz formula we have

$$h^0(C,E\omega)_+-h^0(C,E\omega)_-=0,$$

that means

$$h^0(C,E\omega)_+=h^0(C,E\omega)_-=1.$$

$\square$

\bigskip

\textit{Continuation of the proof of Theorem \ref{yo}:}

Thanks to Lemmas \ref{shit} and \ref{anto} we can extend the morphism $j$ to the closure
$C_E=\overline{\varphi^{-1}(E)}$. We send every point $p\in D(E)$ on the section $s_p$ that vanishes in $p$. We
will denote this morphism

$$\tilde{j}:C_E \lra \pr^1.$$

The conic $C_E$ is either smooth, or the union of two disjoint lines, or a double line. The morphism $\tilde{j}$
is surjective on $\pr^1$ and its degree is 1. Then it exists a morphism $\tau$ s.t. $ \tilde{j}\circ \tau =
Id_{\pr^1}$. This means that we have an isomorphism between $\pr^1=\pr H^0(C,E\omega)$ and a component of $C_E$.
The morphism $\tau$ is equivariant under the action of $\lambda$ and by Lemma \ref{ss} the space
$H^0(C,E\omega)$ has two eigen-spaces. This means that the component that is image of $\tau$ must cut $\pru$ in
two different points, that means that $C_E$ is a smooth conic. $\square$

\begin{thm}\label{hl2}
Let $E$ be a strictly semi-stable vector bundle, then $C_E$ is singular. If $E\cong L\oplus L^{-1}$ for $L\in
JC[2]/\mathcal{O}_C$ then $C_E$ is a double line.
\end{thm}

\textit{Proof:} If the bundle $E$ is strictly semi-stable we know, by Lemma \ref{lem:ber} and Remark
\ref{re:marco} that the fiber consists of two lines so it is either a rank 2 or a rank 1 conic. Moreover the
fibers over the points of $JC[2]/\mathcal{O}_C$ are the double lines $\overline{w_iw_j}$, for $w_i,w_j$ two
different Weierstrass points. These are all the $\lambda$-invariant couples of points, that means that the 15
2-torsion points are the rank 1 locus. $\square$\\

\section{The conic bundle}

In the following, we will often denote $\pr^2_{\omega}$ the linear system $|\omega^2|$.
We will now define a rank 2 projective bundle $\mathcal{E}$  on $\prd$ strictly connected to the classifying map.\\


Since the rational map

$$\Delta: \su \dashrightarrow \prd$$

is surjective, every point of $\prd$ can be represented by a divisor $\Delta(E)$ for a semi-stable bundle $E$ on
$C$. We start by constructing a rank 2 vector bundle $\mathcal{A}$ on $\prd$. Let us first define the fiber
$\mathcal{A}_E$ over the point $\Delta(E)$: we want that $\mathcal{A}_E\subset H^0(C,\omega^3)^*_+$ and that its
dual is the cokernel of the natural multiplication map

\begin{equation}\label{wl}
0 \lra H^0(C,\omega) \stackrel{+\Delta(E)}{\lra} H^0(C,\omega^3)_+ \lra \mathcal{A}_E^*\lra 0.
\end{equation}

Furthermore we can generalize the sequence \ref{wl} to an exact sequence (in fact a global version of the one
just defined) of vector bundles on $\prd$. In order to do this we define a new rank 2 vector bundle
$\mathcal{G}$ on $\pr^2$.

We remark in fact that it exists a natural morphism of vector bundles

$$\nu:\OO_{\prd}(-1) \lra H^0(C,\omega^2)\otimes \OO_{\prd}$$

that sends the fiber over one point $p$ on the line in $H^0(C,\omega^2)$ whose projectivized is $p$ . Now we
twist by $H^0(C,\omega)$ the morphism $\nu$: we get the morphism

$$ Id_{H^0(C,\omega)}\otimes \nu=:\nu^{\prime} :H^0(C,\omega) \otimes \OO_{\prd}(-1) \lra H^0(C,\omega) \otimes H^0(C,\omega^2)\otimes \OO_{\prd}.$$

Moreover it exists a natural multiplication morphism

$$\mu:H^0(C,\omega) \otimes H^0(C,\omega^2)\otimes \OO_{\prd} \lra  H^0(C,\omega^3)_+\otimes \OO_{\prd}.$$

We define

$$\mathcal{G}:=H^0(C,\omega) \otimes \OO_{\prd}(-1)$$

and the injective morphism

$$\alpha:= \mu \circ \nu^{\prime}: \mathcal{G} \lra H^0(C,\omega^3)_+\otimes \OO_{\prd}.$$

The bundle $\mathcal{G}$ is in fact a sub-bundle of $\OO_{\prd}\otimes H^0(C,\omega^3)_+$: the fiber over every
$\Delta(E)\in \prd$ is composed by the divisors of $H^0(C,\omega^3)_+$ of the form $\Delta (E) + \delta$, with
$\delta \in H^0(C,\omega)$. We define

$$\mathcal{A}^*:=coker (\alpha: \mathcal{G} \lra H^0(C,\omega^3)_+\otimes \OO_{\prd})$$

and we get the following exact sequence of bundles on $\prd$.

\begin{equation}\label{ff}
0\lra \mathcal{G} \lra \OO_{\prd}\otimes H^0(C,\omega^3)_+ \lra \mathcal{A^*}\lra 0.
\end{equation}

\begin{defi}
We define the vector bundle $\mathcal{E}$ on $\prd$ as

$$\mathcal{E}:=\mathcal{A}\oplus \OO_{\prd}.$$

\end{defi}

The situation is resumed in the following diagram.

$$\begin{array}{ccc}
\pr \mathcal{E} & \hookrightarrow & |\omega^3|^* \times \prd \\
&&\\
\stackrel{|}{\downarrow} & &\stackrel{|}{\downarrow} \\
&&\\
\pr \mathcal{A} & \hookrightarrow & \pr Sym^3 V^* \times \prd \\
&&\\
 & \searrow & \downarrow \\
 &&\\
 &  & \prd \\
\end{array}$$

Let $y \in \prd$. Then we can identify $y$ and a divisor $a+b$ of degree 2 on $X$, with $a,b\in X$. The fiber
$\pr \mathcal{A}_y$ is the secant line $\overline{ab}$ to $X\subset \pr Sym^3 V^*$ and $\pr \mathcal{E}_y$ is
the plane $<\overline{ab},x> \subset |\omega^3|^*$.\\

Let

$$pr_x: \pro \dashrightarrow \pru$$

be the projection with centre $x$.

\begin{prop}

We have a commutative diagram of rational maps

$$\begin{array}{ccc}
\pro & \stackrel{\varphi}{\dashrightarrow} & |\mathcal{I}_C(2)|^*=\pr^3 \\
\stackrel{|}{\downarrow} pr_x & & \stackrel{|}{\downarrow} \Delta  \\
\pru & \stackrel{\phi}{\dashrightarrow} & \prd\cong |\mathcal{I}_X(2)|^*,
\end{array}$$

where $\phi$ is defined as follows. Given a $t\in \pr^3 \setminus X$, let $l_t$ be the only secant line to $X$
passing by $t$. The application $\phi$ sends $t$ on the pencil $\pr^1\subset|\mathcal{I}_X(2)|$ given by the
quadrics vanishing on the union $X\cup l_t$.
\end{prop}

\textit{Proof:} First we show that there exists a unique secant line $l$ to $X$ that passes by $t$. Projecting
$X$ with centre $t$ we remark that the image is a plane cubic, that has one knot by the genus formula. This
implies that there exists a unique secant line to $X$ passing by $t$. Moreover the projection $pr_x$ induces an
isomorphism

$$pr_x^*|\mathcal{I}_X(2)|\cong |I_S(2)|.$$

Since $\varphi^{-1}(\OO)=S$ (Prop. \ref{pr:giu}) and $\Delta$ is the projection with centre $[\OO\oplus \OO]$,
the diagram commutes. $\square$

\begin{rem}
The map $\phi$ can be defined in a different way. Since\\ $\prd \cong \pr Sym^2 H^0(C, \omega)$, $\phi$ is the
map that sends $t\in \pru \setminus X$ on the pair of points of $X$ cut out by the unique secant line $l_t$ to
$X$ passing by $t$.
\end{rem}




Let $Bl_X\pru$ be the blow-up of $\pru$ along the twisted cubic and

$$\mu:Bl_X\pru \lra \pru$$

the projection on $\pru$. Since $X$ is scheme-theoretically defined by the 3-dimensional space of quadrics
$\mathcal{I}_X(2)$ it exists a morphism $\tilde{\phi}$ that makes the following diagram commute.

$$\begin{array}{ccc}
Bl_X\pru & & \\
\mu\downarrow & \stackrel{\tilde{\phi}}{\searrow}  & \\
X\subset\pru & \stackrel{\phi}{\dashrightarrow} & \prd \\
\end{array}$$

Hence the morphism $\tilde{\phi}$ defines a $\pr^1$-fibration  on $\prd$. Futhermore the exceptional divisor
$E\subset Bl_X\pru$ is the projective bundle $\pr (N_{X|\pru})$ of the normal bundle of $X\subset\pru$.

\begin{lem}
We have an isomorphism

$$N_{X|\pr^3}\cong \OO_{\pr^1}(5) \oplus \OO_{\pr^1}(5).$$
\end{lem}

\textit{Proof:} Let

\begin{eqnarray*}
i:\pr V & \lra & \pru;\\
{[u:v]} & \mapsto & [u^3:u^2v:vu^2:v^3];
\end{eqnarray*}

be the Veronese embedding. We have the following exact sequence.

\begin{equation}\label{edy}
0\lra T_X \lra i^*T_{\pru} \lra N_{X|\pr^3}\lra 0.
\end{equation}

Since $X\cong\pr^1$ we have $T_X\cong \OO_{\pr^1}(2)$. Then we pull-back via $i^*$ the Euler exact sequence and
we get

$$0{\lra} \OO_{\pr^1} \stackrel{k}\lra \OO_{\pr^1}(3)^{\oplus 4}\stackrel{h}{\lra} i^* T_{\pru}\lra 0.$$

Let $l$ be a local section of $\OO_{\pr^1}$ and $u,v$ the coordinates on $\pr^1$, then we can write the morphism
$k$ down in the following way

\begin{eqnarray*}
k:\OO_{\pr^1} & \lra & \OO_{\pr^1}(3)^{\oplus 4};\\
l & \mapsto & (u^3l,u^2vl,uv^2l,v^3l).
\end{eqnarray*}

We denote $X,Y,Z,T$ the coordinates on $\OO_{\pr^1}(3)^{\oplus 4}$. Moreover the morphism $h$ is given by the
equations of the line image of $\OO_{\pr^1}$ in $\OO_{\pr^1}(3)^{\oplus 4}$. Therefore we have

\begin{eqnarray*}
h:\OO_{\pr^1}(3)^{\oplus 4} & \lra & i^*T_{\pru};\\
(X,Y,Z,T) & \mapsto & (vX-uY,vY-uZ,vZ-uT).
\end{eqnarray*}

Hence we have $i^*T_{\pru}\cong \OO_{\pr^1}(4)^{\oplus 3}$ and we can rewrite the exact sequence (\ref{edy}) in
the following way

$$0 \lra \OO_{\pr^1} (2)\cong T_X \stackrel{di}\lra \OO_{\pr^1}(4)^{\oplus 3}\cong i^*T_{\pru} \lra N_{X|\pr^3}\lra 0,$$

where $di$ is the differential of $i$. On the affine open set $\{v\neq 0\}$ the morphism $di$ is defined by the
following equations

\begin{eqnarray*}
di:\OO_{\pr^1} (2) & \lra & \OO_{\pr^1}(4)^{\oplus 3};\\
l & \mapsto & (3lu^2,2luv,lv).
\end{eqnarray*}

Let $(C,D,F)$ be the coordinates on $\OO_{\pr^1}(4)^{\oplus 3}$, then the equations of the line image of
$\OO_{\pr^1}(2)$ in $\OO_{\pr^1}(4)^{\oplus 3}$ are

$$(vA-uB,vB-uC);$$

this implies that $N_{X|\pr^3}\cong \OO_{\pr^1}(5) \oplus \OO_{\pr^1}(5).\square$

\bigskip

Since $\OO_{\pr^1}(5) \oplus \OO_{\pr^1}(5)\cong (\OO_{\pr^1} \oplus \OO_{\pr^1}) \otimes \OO_{\pr^1}(5)$, we
have

$$\pr(N_{X|\pr^3})\cong X \times \pr^1 \cong \pr^1 \times \pr^1.$$

We will denote

$$\sigma: \pr \mathcal{A} \lra \prd$$

the projection of the projective bundle and $E$ the exceptional divisor in $Bl_X\pru$.

\begin{prop}\label{pr:a}

There exists an isomorphism of projective bundles on $\prd$

$$\begin{array}{ccc}
\pr \mathcal{A} & \stackrel{\sim}{\lra} &  Bl_X\pru\\
\sigma \searrow & & \swarrow_{\tilde{\phi}} \\
 & \prd & \\
\end{array}$$

\end{prop}

\textit{Proof:} Let $x\in \prd$, $a+b$ the divisor on $X$ that corresponds to $x$ and $t\in \pr \mathcal{A}_x$.
We recall that $\mathcal{A}$ is a sub-bundle of $\OO_{\pru}\otimes H^0(\omega^3)_+^*$ and the projectivization
gives us an embedding

\begin{eqnarray*}
k:\pr \mathcal{A} & \hookrightarrow & \prd \times \pru.\\
\end{eqnarray*}

Moreover we have

$$k(\pr \mathcal{A}_x)=\overline{ab}\subset \pru.$$

We will often consider $\pr \mathcal{A}$ as a sub-variety of $\prd \times \pru$. If $t\in \pru \setminus X$ then
$\overline{ab}$ is the only secant line to $X$ passing by $t$.\\

We define a morphism

\begin{eqnarray*}
\varpi:= \tilde{\phi} \times \mu : Bl_X\pru & \lra & \prd \times \pru.\\
\end{eqnarray*}

The morphism $\varpi$ has a birational inverse

$$\varpi^{-1}:\pr\mathcal{A}\subset \prd \times \pru \dashrightarrow Bl_X\pru$$

defined as follows. We define $\varpi^{-1}$ on the open set of $\pr\mathcal{A}$ given by the couples $(x,t)$
s.t. $x=a+b$ is a point of $\prd$ s.t. $a\neq b$ and $t\in \{\pr\mathcal{A}_x\setminus X\}$. We send
$(x,t)\in\pr \mathcal{A}$ on $\mu^{-1}(t)\in Bl_X\pru$. Then by Zariski's main theorem the morphism
$\varpi$ induces an isomorphism between $Bl_X\pru$ and $\pr\mathcal{A}$. $\square$\\

Let

$$\beta: \pr (\mathcal{A}\oplus \OO_{\prd}) \dashrightarrow \pr \mathcal{A}$$

be the natural projection and $\eta$ the composed map

$$\eta: \pr \mathcal{E} \stackrel{\beta}{\dashrightarrow} \pr \mathcal{A} \stackrel{\sigma}{\lra} \prd.$$

Since the map $\varphi$ is given by the quadrics in the ideal of $C\subset \pro$, there exists a morphism

$$\overline{\varphi}:Bl_C\pro \lra |2\Theta|$$

that makes the following diagram commute.

$$\begin{array}{ccc}
Bl_C\pro  &  &  \\
&&\\
\downarrow & \searrow^{\overline{\varphi}} & \\
&&\\
\pro & \stackrel{\varphi}{\dashrightarrow} & |2\Theta|\\
\end{array}$$

We will denote $\pr^3_{\OO}$ the blow-up of $|2\Theta|=\pr^3$ at the point $[\OO\oplus \OO]$ and $pr_0 $ the
morphism that resolves the projection $\Delta$ with centre $[\OO\oplus \OO]$ and that makes the following
diagram commute.

$$\begin{array}{ccc}
\pr^3_{\OO}  &  &  \\
&&\\
\downarrow & \searrow^{pr_0} & \\
&&\\
\pr^3 & \stackrel{\Delta}{\dashrightarrow} & \prd\\
\end{array}$$

The morphism

$$pr_0:\pr_{\OO}^3\lra \prd$$

defines a $\pr^1$-fibration  on $\prd$ hence $\pr_{\OO}^3\cong \pr M$ for some rank 2 vector bundle $M$ on
$\prd$. We denote $F= \pr \mathrm{T}_{\OO}\pr^3 \cong \prd$ the exceptional divisor over the origin in
$\pr^3_{\OO}$ . The vector bundle $M$ is defined up to a line bundle $L$, because, as projective varieties $\pr
M \cong \pr (M\otimes L)$, so we choose $M$ once and for all as the vector bundle s.t.

$$\OO_{\pr M}(1) = \OO_{\pr_{\OO}^3}(F).$$

\begin{lem}

We have the equality

$$M=\OO_{\prd} \oplus \OO_{\prd}(1).$$

\end{lem}

\textit{Proof:} We have $M^*=pr_{0*}\OO_{\pr^3_{\OO}}(F)$ and we consider the restriction exact sequence

\begin{equation}\label{ja}
0\lra \OO_{\pr^3_{\OO}}\lra \OO_{\pr^3_{\OO}}(F) \lra \OO_F(F) \lra 0.
\end{equation}

We push down via $pr_{0}$ the exact sequence (\ref{ja}) and we get

$$0\lra \OO_{\prd} \lra M^* \lra \OO_{\prd}(-1)\lra 0.$$

This means that $M^*$ determines an extension class $(e)$ in $Ext^1(\OO(-1),\OO)$. We remark that

$$Ext^1(\OO_{\prd}(-1),\OO_{\prd})=H^1(\prd,\OO(1))=\{0\}$$

thus $M^*$ is the trivial extension, i.e.

$$M^*=\OO_{\prd} \oplus \OO_{\prd} (-1).$$

$\square$

\bigskip

By Proposition \ref{pr:giu}, we have $\overline{\varphi^{-1}(\OO)}=S$, so there exists a morphism

$$\tilde{\varphi}: Bl_S \pro \lra \pr^3_{\OO}$$

that makes the following diagram commute.

$$\begin{array}{ccc}
 Bl_S\pro & & \\
 &&\\
\downarrow & \searrow^{\tilde{\varphi}} &\\
&&\\
Bl_C \pro  &  & \pr^3_{\OO} \\
&&\\
\downarrow & \searrow^{\overline{\varphi}} & \downarrow  \\
&&\\
\pro & \stackrel{\varphi}{\dashrightarrow} & |2\Theta|\\
\end{array}$$

We will denote $\varrho$ the composed map

$$\varrho:Bl_S\pro \stackrel{\tilde{\varphi}}{\lra} \pr^3_{\OO} \stackrel{pr_0}{\lra} \prd$$

and $\pi$ the projection

$$\pi:Bl_S\pro \lra \pro.$$

\begin{prop}

Let $S\subset \pro$ be the cone over the twisted cubic $X$ of Proposition \ref{pr:giu}. Then there exists an
isomorphism of projective bundles on $\prd$

$$\begin{array}{ccc}
\pr \mathcal{E} & \stackrel{\sim}{\lra} &  Bl_S\pro\\
&&\\
\eta \searrow & & \swarrow\varrho\\
&&\\
 & \prd. & \\
\end{array}$$

\end{prop}
\textit{Proof:} Let $x\in \prd$, $a+b$ the divisor on $X$ corresponding to $x$ and $s\in \pr \mathcal{E}_x$. We
recall that $\mathcal{E}$ is a sub-bundle of $\OO_{\pru}\otimes H^0(\omega^3)^*$, thus we have an embedding

\begin{eqnarray*}
j:\pr\mathcal{E} & \hookrightarrow & \prd \times \pro.\\
\end{eqnarray*}

Moreover we have

$$j(\pr \mathcal{E}_x)=x \times \langle a + b + \lambda(a) + \lambda (b) \rangle \subset \pro.$$

We will often consider $\pr \mathcal{E}$ as a sub-variety of $\prd \times \pro$. We also remark that

$$j(\pr \mathcal{E}_x)\cap \prd \times \pru= k(\pr \mathcal{A}_x)=x \times \overline{ab}.$$

We define a morphism

\begin{eqnarray*}
\varpi^{\prime}:=(\varrho,\pi): Bl_S\pro  & \lra & \prd \times \pro .\\
\end{eqnarray*}

The morphism $\varpi^{\prime}$ has a birational inverse

$$\varpi^{\prime -1}: \pr\mathcal{E}\subset \prd \times \pro \dashrightarrow Bl_S\pro$$
 defined as follows. We define $\varpi^{\prime -1}$ on the open set of $\pr\mathcal{E}$ given by the couples
$(x,s)$ s.t. $x=a+b$ is a point of $\prd$ s.t. $a\neq b$ and $s\in \{\pr\mathcal{E}_x\setminus S\}$. The pair
$(x,s)$ is sent on $\pi^{-1}(s)\in Bl_S\pro$. Then, by Zariski's main theorem $\varpi^{\prime}$ induces an
isomorphism between $Bl_S\pro$ and
$\pr\mathcal{E}.\square$\\

We recall that we denoted $E$ the exceptional divisor of $Bl_X\pru$.

\begin{thm}\label{hate}
The restricted map

$$\tilde{\phi}_{|E}:E \lra \prd$$

is a morphism of degree 2 ramified along the conic that is the image of $\pr V$ via the quadratic Veronese
embedding $Ver_2$.
\end{thm}

\textit{Proof:} We recall that $E\cong X\times \pr^1$. Moreover we remark that $\tilde{\phi}_{|E}$ is given by
the differential of $\varphi$ and we have

$$\tilde{\phi}_{|{\{a\}}\times \pr^1}: \pr (N_{X|\pr^3,a}) {\lra} \prd.$$

Let $a,b\in X$, then we have

$$\tilde{\phi} ( \overline{ab}-\{a,b\})= a+b \in |\omega^2|.$$

Let $v_b\in \mathrm{T}_a\pr^3$ the tangent vector to $\pr^3$ with direction $\overline{ab}$. Then, since the
line $\overline{ab}$ is contracted to a point, we have

$$\tilde{\phi}(v_b)=a+b \in |\omega^2|.$$

Furthermore, every normal vector $v\in \pr^1 \cong \pr (N_{X|\pr^3,a})$ is of type $v_b$ for a point $b\in X$.
This implies that

$$\tilde{\phi}(\{a\} \times \pr^1) =D_a:=\{a+b|b\in X\}.$$

The line $D_a\subset \prd$ is the tangent line at the point $2a$ to the conic in $\prd$ obtained as the image of
the Veronese morphism

\begin{eqnarray}
Ver_2: \pr V & \lra & |\OO_{\pr^1}(2)|^*=\prd;\\
p & \mapsto & 2p.
\end{eqnarray}

Let $a+b$ be again the divisor on $X$ corresponding to $x \in \prd$, then the fiber $\tilde{\phi}^{-1}(x)$ in $X
\times \pr^1\cong \pr (N_{X|\pr^3})$ is composed by two points $\{(a,\alpha),(b,\beta)\}$ if $x$ is not
contained in the conic. We also remark that if $x$ is a point contained in the conic the fiber is just one
point. This defines a degree
2 covering ramified along the conic.$\square$\\

\bigskip

We will denote $\tilde{E}$ the exceptional divisor of $Bl_{S}\pro \cong \pr \mathcal{E}$. We have

$$\tilde{E}\cong\overline{\beta^{-1}(E)}.$$

\begin{prop}\label{hl}
The restricted map

$$\tilde{\varphi}_{|\tilde{E}}:\tilde{E} \lra E \stackrel{2:1}{\lra} F \cong \prd $$

defines a conic bundle
\end{prop}

\textit{Proof:} The situation is the following

$$\begin{array}{cccc}
\tilde{E} \subset & \pr (\mathcal{A} \oplus \OO) = Bl_{S}\pro & &\\
\stackrel{|}{\downarrow}\beta  &  &\searrow \tilde{\varphi} &\\
 &&&\\
 E \subset & \pr \mathcal{A} = Bl_X \pru & \longrightarrow & \pr^3_{\OO}\\
 &&&\\
\sigma \downarrow  & & &\\
 &&&\\
\prd &  & & \\
\end{array}$$

Let us consider the composed map

$$\eta_{|\pr \tilde{E}}:= \sigma \circ \beta_{|\pr \tilde{E}}: \pr \tilde{E} \lra \pr^2.$$

The fiber of $\eta$ over a point $x\in \pr^2$ is a rank two conic if $x\not \in Ver_2(\pr^1)$ and a double line
if
$x\in Ver_2(\pr^1). \square$\\

Proposition \ref{hl} and Theorem \ref{hl2} imply the following theorem.

\begin{thm}\label{yooo}
The morphism

$$\tilde{\varphi}:Bl_S(\pr^4)\lra \pr^3_{\OO}$$

is a conic bundle whose discriminant locus is the blow-up at the origin of the Kummer surface $K^0$.
\end{thm}

\begin{rem}
Moreover we remark that the conic $Ver_2(\pr^1)$ is the tangent cone at the origin of the Kummer surface.
\end{rem}

\begin{rem}
We recall that $Pic(\pr \mathcal{A}) \cong \mathbb{Z}^2$, notably

$$Pic(\pr \mathcal{A}) = \OO_{\pr\mathcal{A}}(1) \mathbb{Z} \times \sigma^*\OO_{\prd}(1) \mathbb{Z}.$$

\end{rem}

We also recall that the map

\begin{eqnarray*}
\mu:\pr \mathcal{A} & \lra & \pru = \pr Sym^3 V;\\
\pr \mathcal{A}_x & \mapsto & \overline{ab};\\
\end{eqnarray*}

that sends the fiber over $x$ on the secant line $\overline{ab}$ to $X$ is the projection of the blow-up of
$\pru$ along $X$, hence $\mu^{-1}(X)=E$.

\begin{lem}

We have

$$\mu^*\OO_{\pru}(1) = \OO_{\pr \mathcal{A}}(1).$$
\end{lem}

\textit{Proof:} Our aim is to determine two integers $l,k\in \mathbb{Z}$ s.t.

$$\mu^*\OO_{\pru}(1)=\OO_{\pr\mathcal{A}}(l)\otimes \sigma^*\OO_{\prd}(k).$$

Since the map $\mu$ is the projection of the blow-up of $\pru$ along $X$, we have $l=1$. Now we need to
determine $k$. We have

$$H^0(\pru,\mu^*\OO_{\pru}(1))=H^0(\pru,\mu_*\mu^*\OO_{\pru}(1)),$$

and by projection formula this is equal to $H^0(\pru,\OO_{\pru}(1)\otimes \mu_*\OO_{\pr\mathcal{A}}).$ Since the
fibers of $\mu$ are connected

$$\mu_*\OO_{\pr\mathcal{A}} = \OO_{\pru}.$$

Therefore

$$H^0(\pr\mathcal{A},\mu^*\OO_{\pru}(1))=H^0(\pru,\OO_{\pru}(1))=Sym^3V.$$

By taking the cohomology of the exact sequence (\ref{ff}) we get an isomorphism

\begin{equation}\label{fly}
H^0(\prd, \mathcal{A}^*) \cong Sym^3 V.
\end{equation}

In order to determine $k$ we compute $H^0(\pr\mathcal{A},\OO_{\pr\mathcal{A}}(1)\otimes \sigma^*\OO_{\prd}(k))$.
By the projection formula we get

\begin{equation}\label{dp}
H^0(\pr\mathcal{A},\OO_{\pr\mathcal{A}}(1)\otimes \sigma^*\OO_{\prd}(k))=H^0(\prd, \mathcal{A}^*\otimes
\OO_{\prd}(k)).
\end{equation}

We twist the exact sequence (\ref{ff}) by the line bundle $\OO_{\prd}(k)$ and we take the cohomology of the
obtained sequence. This leads us to conclude that the equality

$$H^0(\prd,\mathcal{A}^*\otimes\OO_{\prd}(k))= Sym^3V$$

is possible only for $k=0. \square$\\

Now we compute the class of the exceptional divisor $E$ in the Picard group of $\pr \mathcal{A}$.

\begin{thm}
We have an isomorphism in $Pic(\pr \mathcal{A})$

$$\OO(E)_{\pr \mathcal{A}} \cong \OO_{\pr \mathcal{A}}(2) \otimes \sigma^* \OO_{\prd}(-1).$$
\end{thm}

\textit{Proof:} First we compute the first factor. We recall that the restricted map

$$\til:E\lra \prd$$

defines a degree 2 morphism ramified along a smooth conic (Thm. \ref{hate}). If $x\in \prd$ then the
intersection $E_x = \tilde{E} \cap \pr \mathcal{A}_x$ is made up of two points. This implies that

$$\OO(E)_{\pr \mathcal{A}} \cong \OO_{\pr \mathcal{A}}(2) \otimes \sigma^* \OO_{\prd}(m),$$

for some $m\in \mathbb{Z}$. We have the following equalities.

$$H^0(\pr\mathcal{A}, \OO(1))= H^0(\prd, \sigma_* \OO_{\pr \mathcal{A}} (1)) = H^0(\prd, \mathcal{A}^*).$$

Let $a,b\in X$ and let $x\in\prd$ be as usual the point corresponding to the divisor $a+b$.  
Let us consider a smooth quadric $Q\subset |\mathcal{I}_X(2)|.$ Then we have

$$\mu^{-1}(Q)=E + R_Q \subset |\OO_{\pr\mathcal{A}}(2)|,$$

where we denote by $R_Q\subset\pr\mathcal{A}$ the residual divisor. Moreover, we have either
$\overline{ab}\subset Q$, or $\overline{ab} \cap Q = \{a,b\}$.

We define

$$\mathcal{C}_Q:=\{a+b=x\in \pr^2|\overline{ab}\subset Q\}$$

and $R_Q=\sigma^{-1}(C_Q)$.

It is well known that, since $Q$ is smooth, we have an isomorphism $Q \cong \pr^1 \times \pr^1$ via the Segre
embedding and that $X \subset Q$ can be seen as the zero locus of a bihomogeneus polynomial of degree $(1,2)$.

This means that we have an embedding

\begin{eqnarray*}
X=\pr^1 & \hookrightarrow & Q = \pr^1 \times \pr^1;\\
{[u:v]} & \mapsto & ([u^2,v^2],[u,v]).\\
\end{eqnarray*}

This implies that, if we choose a $p\in \pr^1$ and we let $t$ vary in the other $\pr^1$, the lines of the ruling
$\{t\}\times p$ intersect $X$ in two points. The lines of the other ruling of $Q$ intersect $X$ in just one
point. Let $\alpha$ be the following morphism

\begin{eqnarray*}
\alpha:\pr^1 & \lra & \pr^1;\\
{[u:v]}& \mapsto & [u^2:v^2].
\end{eqnarray*}

Then we have a linear embedding

$$H^0(X,\OO(1)) \stackrel{\alpha^*}{\lra} H^0(X,\OO(2))$$

and the line $\pr(\alpha^*(H^0(X,\OO(1))))$ is $\mathcal{C}_Q$. This implies that

$$\OO_{\pr\mathcal{A}} (R)= \sigma^*\OO_{\prd}(1).$$

and thus

$$\OO_{\pr\mathcal{A}}(E)=\OO_{\pr\mathcal{A}}(2)\otimes \sigma^*\OO_{\prd}(-1).$$
$\square$

The morphism $\tilde{\varphi}$ then makes the following diagram of projective bundles on $\pr^2$ commute.

$$\begin{array}{ccc}
Bl_S \pro \cong \pr (\mathcal{A} \oplus \OO_{\prd}) & \stackrel{\tilde{\varphi}}{\lra} & \pr [\OO_{\prd} \oplus \OO_{\prd}(1)]=\pr^3_{\OO}\\
&&\\
\varrho\searrow & & \swarrow pr_0\\
&&\\
 & \prd & \\
\end{array}$$

The pull back $\tilde{\varphi}^*$ induces a homomorphism

\begin{equation}\label{wft}
\tilde{\varphi}^*: Pic(\pr^3_{\OO}) \lra Pic(\pr (\mathcal{A}\oplus \OO_{\prd})),
\end{equation}

and both Picard groups are isomorphic to $\mathbb{Z}\times \mathbb{Z}$.

\begin{prop}
The homomorphism of equation (\ref{wft}) is the following.
\begin{eqnarray*}
\mathbb{Z}\times \mathbb{Z} & \lra & \mathbb{Z}\times \mathbb{Z}\\
(a,b) & \mapsto & (2a,b-a).
\end{eqnarray*}
\end{prop}

\textit{Proof:} We have the following equalities:

\begin{eqnarray*}
Pic(\pr^3_{\OO})= \ZZ \OO_{\pr^3_{\OO}}(1) \times \ZZ pr_0^* \OO_{\prd}(1);\\
Pic(\pr (\mathcal{A}\oplus \OO))= \ZZ \OO_{\pr (\mathcal{A}\oplus \OO)}(1) \times \ZZ \eta^* \OO_{\prd}(1).\\
\end{eqnarray*}

Moreover, since $pr_0\circ \tilde{\varphi} = \eta$, we have

$$\tilde{\varphi}(0,b) = (0,b).$$

We recall that we chose $F$ in such a way that $\OO_{\pr^3_{\OO}}(1)=\OO_{\pr^3_{\OO}}(F)$, this means that we
have

$$\tilde{\varphi}^* \OO_{\pr^3_{\OO}}(1) = \tilde{\varphi}^* \OO_{\pr^3_{\OO}}(E)=\OO_{\pr\mathcal{E}}(\tilde{\varphi}^{-1}(F)) = \OO_{\pr\mathcal{E}}(\beta^{-1}(E))= \OO_{\pr\mathcal{E}}(2) \otimes \eta^* \OO_{\prd}(-1).$$

This implies

$$(a,0)\mapsto (2a,-a).$$
$\square$\\

As a consequence of the last proposition, we have

\begin{equation}\label{side}
\tilde{\varphi}^* \OO_{\pr^3_{\OO}}(1) = \OO_{\pr\mathcal{E}}(2) \otimes \eta^*\OO_{\prd}(-1).
\end{equation}

Furthermore $\tilde{\varphi}$ induces a natural morphism

\begin{equation}\label{stro}
\OO_{\pr^3_{\OO}}(1) \lra \tilde{\varphi}_* \tilde{\varphi}^* \OO_{\pr^3_{\OO}}(1).
\end{equation}

By applying $pr_{0*}$ to the morphism (\ref{stro}) and using the equality (\ref{side}), we obtain a morphism of
sheaves on $\prd$

$$\tilde{\varphi}^*: pr_{0*}\OO_{\pr^3_{\OO}}(1)= \OO \oplus \OO(-1)=M^* \lra \eta_*(\OO_{\pr\mathcal{E}}(2) \otimes \eta^* \OO_{\prd}(-1)).$$

By the projection formula,

$$\eta_*(\OO_{\pr\mathcal{E}}(2) \otimes \eta^* \OO_{\prd}(-1))\cong Sym^2 (\mathcal{A}^*\oplus \OO_{\prd}) \otimes \OO_{\prd}(-1).$$

\begin{rem}
We showed that $\pr M^*$ defines a $\pr^1$-bundle on $\prd$. Let $y\in \prd$ and let $c,d$ be the points of $X$
s.t. $y$ is the divisor $c+d$ on $X$. We recall that $C$ is a degree 2 covering of $X\cong \pr V$. Then $\pr
(\mathcal{A}\oplus \OO)_y\subset \pro$ is the $\pr^2$ generated by the two pairs of points of $C$ whose images
in $X$ are respectively $c$ and $d$. The fiber $\pr M^*_y$ is in fact the pencil of conics in $\pr
\mathcal{E}_y$ that pass by these four points.
\end{rem}

\begin{lem}\label{db}

We have

\begin{eqnarray*}
h^0(\prd, \mathcal{A}^*(-1))=0;\\
h^1(\prd, \mathcal{A}^*(-1))=0.
\end{eqnarray*}

\end{lem}

\textit{Proof:} We twist the exact sequence (\ref{ff}) by the vector bundle $\OO_{\prd}(-1)$ and we find

$$0 \lra V \otimes \OO_{\prd}(-2) \lra \OO_{\prd}(-1) \otimes Sym^3V \lra \mathcal{A}^*(-1)\lra 0. $$

By taking cohomology, we get

$$h^0(\prd, \OO(-1))=h^1(\prd, \OO(-2))=0.$$

This gives us the first equality. Then, we have

$$h^1(\prd, \OO(-1))=0$$

and by duality

$$H^2(\prd,\OO(-2))\cong H^0(\prd,\OO(-1),$$

hence $h^2(\prd,\OO(-2))=0$. This implies the second equality.$\square$

\begin{lem}

We have

$$h^0(\prd,Sym^2 \mathcal{A}^*)=10.$$

\end{lem}

\textit{Proof:} By twisting the exact sequence (\ref{ff}) by $\mathcal{A}^*$ we obtain the following exact
sequence.

$$0\lra \mathcal{A}^*(-1) \otimes V \lra \mathcal{A}^* \otimes Sym^3 V \lra \mathcal{A}^{*2} \lra 0.$$

This implies that we have

$$H^0(\prd, \mathcal{A}^{*2})\cong H^0( \prd,\mathcal{A}^*)\otimes Sym^3V)= Sym^3V\otimes Sym^3V,$$

and therefore $h^0(\prd, \mathcal{A}^{*2})=16$. Moreover we remark that

$$h^0(\prd, \mathcal{A}^{*2})= h^0(\prd, \wedge^2 \mathcal{A}^*) \oplus h^0(\prd,Sym^2 \mathcal{A}^*).$$

By taking determinants in the exact sequence (\ref{ff}) we get that $\bigwedge^2 \mathcal{A}^*=\OO_{\prd}(2)$
and  $h^0(\prd, \bigwedge^2 \mathcal{A}^*)=h^0(\prd, \OO(2))=6$. This implies directly the lemma.$\square$

\begin{lem}
We have

$$h^0(\prd,Sym^2A^*(-1))=1.$$

\end{lem}

\textit{Proof:} We twist the exact sequence (\ref{ff}) by the line bundle $\mathcal{A}^*(-1)$ and we get the
following exact sequence.

$$0\lra \mathcal{A}^*(-2) \otimes V \lra \mathcal{A}^*(-1) \otimes Sym^3 V \lra \mathcal{A}^{*2}(-1) \lra 0.$$

By passing to cohomology we have that

$$h^1(\prd,\mathcal{A}^*(-1))\otimes Sym^3V=h^0(\prd,\mathcal{A}^*(-1))\otimes Sym^3V=0.$$

Thanks to Lemma \ref{db} we have then

\begin{equation}\label{dak}
H^0(\prd, \mathcal{A}^{*2}(-1))\cong H^1(\prd,\mathcal{A}^*(-2))\otimes V.
\end{equation}

We twist the exact sequence (\ref{ff}) by the line bundle $\OO_{\prd}(-2)$ and we get

$$0 \lra V \otimes \OO_{\prd}(-3) \lra \OO_{\prd}(-2) \otimes Sym^3V \lra \mathcal{A}^*(-2)\lra 0. $$

Moreover we remark that $h^1(\prd,\OO(-2)=0$ and $h^2(\prd,\OO(-3))=1$. By duality
$h^2(\prd,\OO(-2))=h^0(\prd,\OO(-1))=0$, then

$$h^1(\prd,\mathcal{A}^*(-2))=h^2(\prd,\OO(-3))\otimes V=2.$$

The equality (\ref{dak}) implies that $h^0(\prd,\mathcal{A}^{*2}(-1)=4$. The vector bundle
$\mathcal{A}^{2*}(-1)$ decomposes as  $ Sym^2 \mathcal{A}^*(-1) \oplus \bigwedge^2 \mathcal{A}^*(-1)$. Since

$$h^0(\prd,\bigwedge^2 \mathcal{A}^*(-1))=h^0(\prd,\OO(1))=3$$

we have that

$$h^0(\prd,Sym^2 \mathcal{A}^*(-1))=4-3=1.$$
$\square$

\begin{cor}
We have

$$dim (\mathrm{Hom}(\OO_{\prd}\oplus\OO_{\prd}(-1),Sym^2 (\mathcal{A}^*\oplus \OO_{\prd})\otimes \OO_{\prd}(-1)))=16.$$

\end{cor}

\textit{Proof:} We will denote

$$B:=\mathrm{Hom}(\OO_{\prd}\oplus\OO_{\prd}(-1),Sym^2 (\mathcal{A}^*\oplus \OO_{\prd})\otimes \OO_{\prd}(-1)).$$

We have

$$B= H^0(\prd,Sym^2\mathcal{A}^*\oplus \OO_{\prd}\oplus \mathcal{A}^* \oplus Sym^2 \mathcal{A}^*(-1) \oplus
\OO_{\prd}(-1) \oplus \mathcal{A}^*(-1)).$$

Summing the dimensions of the direct summands we have

$$dim B=10+1+4+1=16.$$
$\square$

\begin{rem}
The projective plane $\prd=\pr Sym^2V$ does not depend on the curve $C$, because $V$ is an abstract vector
space. Since it is defined by the exact sequence $\ref{ff}$, even $\mathcal{A}$ does not depend on the curve
$C$.
\end{rem}

This means that to a given a genus 2 curve $C$ we can associate a conic bundle over $\pr^2$, notably the bundle
defined by the section

$$\tilde{\varphi}_C\in \pr B.$$

Its discriminant locus is the blow-up (at the origin) of the Kummmer surface  $K^0=Jac(C)/\pm Id$. In this way
we build a moduli map

\begin{eqnarray*}
\Xi:\{ \mathrm{smooth\ genus\ 2\ curves} \} & \lra & \pr^{15}=\pr B;\\
C & \mapsto & \tilde{\varphi}_C.
\end{eqnarray*}

\begin{lem}\label{klu}
Let $Y$ a smooth projective variety,

$$f:G\lra Y$$

a conic bundle, then there exists a rank 3 vector bundle $F$ on Y, a line bundle $L$ and one section $q\in
H^0(Y,Sym^2 H\otimes L^k)$ for some integer $k$, s.t. G is the zero scheme of $q$ in the projective bundle $\pr
F$.
\end{lem}

\textit{Proof:} The proof follows that of Proposition 1.2 of \cite{bo:ji}. The assertion is equivalent to the
 existence of a line bundle $N$ on $G$ inducing the sheaf $\OO_{G_s}(1)$ on every fiber $G_s$. More precisely we
 will have $H=f_*N$. If we consider the sheaf of differentials of maximum degree $\omega_G$ then by the
 adjunction formula we have
$\omega_{G|G_s}\cong \omega_{G_s}$ and $\omega_{G_s}\cong \OO_{G_s}(-1)$. Then we take $N
=\omega_G^{-1}.\square$

\begin{prop}

Let

$$\tilde{\varphi}:Bl_S\pro \lra \pr^3_{\OO}$$

be the conic bundle of Theorem \ref{yooo}. Then $Bl_S\pro$ is a divisor in the total space of the bundle $\pr
(pr_0^*\mathcal{E})$ on $\pr^3_{\OO}$.
\end{prop}

\textit{Proof:} Let $y\in\pr^3_{\OO}$, then we have

$$\tilde{\varphi}^{-1}(y)\subset \pr(\mathcal{E}_{pr_0(y)}).$$

Then the rank 3 vector bundle associated to $\tilde{\varphi}$ is $pr_0^*\mathcal{E}.\square$

\bigskip

Moreover $Bl_S\pro = \pr \mathcal{E}$. This means we have the following diagram, where the lower square is a
fiber product.

\bigskip

\ \ \ \ \ \ \ \ \ \ \ \ \ \ \ \ \ \ \ \ \ \ \ \ \ \ \ \ \ \ \ 
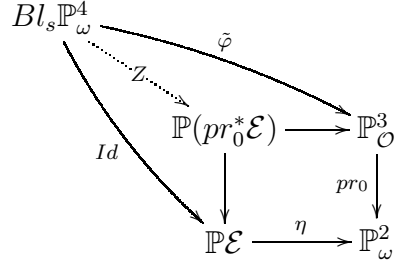
\xymatrix{
  Bl_s\pro \ar@/_/[ddr]_{Id} \ar@/^/[drr]^{\tilde{\varphi}}
    \ar@{.>}[dr]|-{Z}                   \\
   & \pr (pr_0^*\mathcal{E}) \ar[d]^{} \ar[r]
                      & \pr^3_{\OO}\ar[d]_{pr_0}    \\
   & \pr \mathcal{E} \ar[r]^{\eta}     & \prd              }

In this diagram $Z$ is the embedding of $Bl_S\pro$ in $\pr (pr_0^*\mathcal{E})$ induced by the universal
property of fiber product.

\section{Stability and deformations of $\mathcal{A}$.}

\subsection{Stability}

In this section we will go through the question of the stability of $\A$ and we will calculate its space of
deformations. Let $H$ be the hyperplane class that generates $Pic(\pr^2)$. If we take the cohomology of the
exact sequence \ref{ff} (and of its dual sequence) and we compute the Chern polynomials we find that

\begin{eqnarray*}
c_1(\A)=-2H,\\
c_2(\A)=3H^2,\\
c_1(\A^*)=2H,\\
c_2(\A^*)=3H^2.
\end{eqnarray*}

Then we have that the slope $\mu(\A^*)=1$.

\begin{rem}
If we twist by $\OO_{\prd}(1)$ the exact sequence \ref{ff} that defines $\A^*$ we get

\begin{equation}\label{stai}
0\lra \OO_{\prd} \otimes H^0(C,\omega) \lra \OO_{\prd}(1) \otimes H^0(C,\omega^3)_+ \lra \A^*(1) \lra 0.
\end{equation}

According to \cite{dk} a bundle on $\pr^n$ that has a linear resolution like $\A^*(1)$ a \textit{Steiner
bundle}.
\end{rem}

\begin{defi}
Let $K$ be a complex projective manifold with $Pic(K)\cong \ZZ$ and $E$ a vector bundle of rank $r$ on $K$. Then
the bundle $E$ on $K$ is called \textit{normalized} if $c_1(E)\in \{-r+1,\dots,-1,0\}$, i.e. if $-1<\mu(E)\leq
0$. We denote by $E_{norm}$ the unique twist of $E$ that is normalized.
\end{defi}

The following criterion for the stability of vector bundles on $K$ is a consequence of the definition.

\begin{prop}\textbf{(Hoppe)}\cite{hl}\label{kit}
Let $V$ be a vector bundle on a projective manifold $K$ with $Pic(K)\cong \ZZ$. If $H^0(X,\wedge^q V_{norm})=0$
for any $1\leq q \leq rkV-1$, then $V$ is stable.
\end{prop}

\begin{thm}
The vector bundle $\A$ on $\prd$ is stable.
\end{thm}

 \textit{Proof:}We remark that $Pic(\prd)\cong \ZZ$ and that $c_1(\A^*(-1))=0$, this means that $\A^*_{norm}=\A^*(-1)$. By
taking the cohomology of the exact sequence \ref{stai} we remark that $H^0(\A^*(-1))=0$. By proposition
\ref{kit} then $\A^*$ is stable and this in turn implies that $\A$ is stable. $\square$

\subsection{Deformations of $\A$}

Before computing directly the dimension of the space of deformations of $\A$, i.e. $dim (Ext^1(\A,\A))$, we need
some technical lemmas.

\begin{lem}\label{black}
We have

\begin{eqnarray*}
h^1(\prd,\A^*)=0,\\
h^2(\prd,\A^*)=0.
\end{eqnarray*}
\end{lem}

\textit{Proof:} We take the cohomology of the exact sequence \ref{ff} and we get

$$0 \ra H^1(\prd,\A^*) \ra H^2(\prd,\OO(-1) \ra H^2(\prd, \OO) \ra H^2(\prd,\A^*) \ra 0.$$

By duality we have $H^2(\prd,\OO(-1)\cong H^0(\prd,\OO(-2))$ and $H^2(\prd, \OO)\cong H^0(\prd,\OO(-3))$ and
both spaces are zero dimensional. This implies our statement. $\square $\\

Furthermore, we recall that $H^0(\prd,\A^*)\cong H^0(\prd,\omega^3)_+$.

\begin{lem}\label{yl}
We have
\begin{eqnarray*}
h^0(\prd,\A^*(1))=10,\\
h^1(\prd,\A^*(1))=0,\\
h^2(\prd,\A^*(1))=0.
\end{eqnarray*}
\end{lem}

\textit{Proof:} We take the cohomology of the exact sequence \ref{stai} and we get

$$0\ra H^0(C,\omega) \ra H^0(\prd,\OO(1))\otimes H^0(C,\omega^3)_+ \ra H^0(\prd,\A^*(1)) \ra H^1(\prd,\OO)\ra$$

$$\ra H^1(\prd,\OO(1)) \ra H^1(\prd,\A^*(1)) \ra H^2(\prd, \OO) \ra H^2(\prd,\OO(1)) \ra H^2(\prd,\A^*(1)) \ra
0.$$

Since $h^1(\prd,\OO)=0$, we have $h^0(\prd,\A^*(1))=10$. Moreover, by duality we have also

$$h^1(\prd,\OO(1))=h^2(\prd, \OO)= h^2(\prd,\OO(1))=0.$$

Hence $h^1(\prd,\A^*(1))= h^2(\prd,\A^*(1))=0. \square$\\

We are ready to state the main result of this section.

\begin{thm}
The space of deformations of $\A$ has dimension equal to 5.
\end{thm}

\textit{Proof:} We twist by $\A^*$ the dual of the exact sequence \ref{ff} and we obtain the following.

\begin{equation}\label{ev}
0 \lra \A^*\otimes \A \lra \A^* \otimes H^0(C,\omega^3)_+^* \lra \A^*(1) \otimes H^0(C,\omega)^* \lra 0.
\end{equation}

By taking cohomology we get the following long exact sequence

$$0 \ra H^0(\prd,\A^*\otimes \A) \lra H^0(\prd, \A^*)\otimes H^0(C,\omega^3)_+^* \lra H^0(\prd,\A^*(1))\otimes
H^0(C,\omega)^* \ra $$
$$\ra H^1(\prd,\A^*\otimes \A) \ra H^1(\prd, \A^*)\otimes H^0(C,\omega^3)_+^* \ra H^1(\prd,\A^*(1))\otimes
H^0(C,\omega)^* \ra  $$
$$\ra H^2(\prd,\A^*\otimes \A) \ra \ra H^2(\prd, \A^*)\otimes H^0(C,\omega^3)_+^* \ra H^2(\prd,\A^*(1))\otimes
H^0(C,\omega)^* \ra 0. $$\

Lemma \ref{yl} and Lemma \ref{black} imply that $h^2(\prd,\A^*\otimes \A)=0$. Moreover we have that $dim
H^0(\prd,\A^*)\otimes H^0(C,\omega^3)_+^*=16$ and $dim H^0(\prd,\A^*(1))\otimes H^0(C,\omega)^*=20$. This in
turn implies that

$$h^1(\prd,\A^*\otimes \A) - h^0(\prd,A^*\otimes \A)=20-16=4.$$

Now $H^1(\prd,\A^*\otimes \A)\cong Ext^1(\A,\A)$ and, since $\A$ is stable

$$dim Hom(\A,\A)=h^0(\prd,\A^*\otimes \A)=1.$$

This means that $h^1(\prd,\A^*\otimes \A)=5.\square$

\bibliographystyle{amsalpha}
\bibliography{bibkumwed}

\providecommand{\bysame}{\leavevmode\hbox to3em{\hrulefill}\thinspace}
\providecommand{\MR}{\relax\ifhmode\unskip\space\fi MR }
\providecommand{\MRhref}[2]{%
  \href{http://www.ams.org/mathscinet-getitem?mr=#1}{#2}
}
\providecommand{\href}[2]{#2}
\begin{thebibliography}{Bea88}

\bibitem[Bea77]{bo:ji}
Arnaud Beauville, \emph{Vari\'et\'es de {P}rym et jacobiennes
  interm\'ediaires}, Ann. Sci. \'Ecole Norm. Sup. (4) \textbf{10} (1977),
  no.~3, 309--391.

\bibitem[Bea88]{bo:fib1}
\bysame, \emph{Fibr\'es de rang {$2$} sur une courbe, fibr\'e d\'eterminant et
  fonctions th\^eta}, Bull. Soc. Math. France \textbf{116} (1988), no.~4,
  431--448.

\bibitem[Ber92]{ab:rk2}
A.~Bertram, \emph{Moduli of rank-{$2$} vector bundles, theta divisors, and the
  geometry of curves in projective space}, J. Differential Geom. \textbf{35}
  (1992), no.~2, 429--469.

\bibitem[Bol07]{bol:wed}
M.~Bolognesi, \emph{On weddle surfaces and their moduli}, Adv. in Geom.
  \textbf{7} (2007), no.~1, 1--99.

\bibitem[DK93]{dk}
I.~Dolgachev and M.~Kapranov, \emph{Arrangements of hyperplanes and vector
  bundles on {$\bold P\sp n$}}, Duke Math. J. \textbf{71} (1993), no.~3,
  633--664.

\bibitem[DN89]{dn:pfv}
J.-M. Drezet and M.~S. Narasimhan, \emph{Groupe de {P}icard des vari\'et\'es de
  modules de fibr\'es semi-stables sur les courbes alg\'ebriques}, Invent.
  Math. \textbf{97} (1989), no.~1, 53--94.

\bibitem[GH78]{gh:pag}
P.~Griffiths and J.~Harris, \emph{Principles of algebraic geometry},
  Wiley-Interscience [John Wiley \& Sons], New York, 1978, Pure and Applied
  Mathematics.

\bibitem[HL97]{hl}
Daniel Huybrechts and Manfred Lehn, \emph{The geometry of moduli spaces of
  sheaves}, Aspects of Mathematics, E31, Friedr. Vieweg \& Sohn, Braunschweig,
  1997.

\bibitem[LN83]{ln85}
H.~Lange and M.~S. Narasimhan, \emph{Maximal subbundles of rank two vector
  bundles on curves}, Math. Ann. \textbf{266} (1983), no.~1, 55--72.

\bibitem[NR69]{rana:cra}
M.~S. Narasimhan and S.~Ramanan, \emph{Vector bundles on curves}, Algebraic
  Geometry (Internat. Colloq., Tata Inst. Fund. Res., Bombay, 1968), Oxford
  Univ. Press, London, 1969, pp.~335--346.

\end{thebibliography}


\bigskip

Michele Bolognesi\\
Dipartimento di Matematica
Università di Pavia\\
Via Ferrata 1\\
27100 Pavia\\
Italy\\
E-mail: michele.bolognesi@unipv.it

\end{document}